# A filtering technique for the matrix power series being near-sparse

Wu Feng‡; Zhu Li†; Zhao Yuelin†; Zhang Kailing†

†Department of Mechanics, Dalian University of Technology, Dalian 116023, P.R.China

‡School of Mathematical Sciences, State Key Laboratory of Structural Analysis of Industrial Equipment, Dalian University of Technology, Dalian 116023, P.R.China (Tel.:+86-13940846142, Corresponding author: vonwu@dlut.edu.cn).

Project supported by the National Natural Science Foundation of China grants (Nos. 11472076 and 51609034), the Science Foundation of Liaoning Province of China (No. 2021-MS-119), the Dalian Youth Science and Technology Star project (No. 2018RQ06), and the Fundamental Research Funds for the Central Universities grant (Nos. DUT20RC(5)009 and DUT20GJ216).

**Abstract:** This work presents a new algorithm for matrix power series which is near-sparse, that is, there are a large number of near-zero elements in it. The proposed algorithm uses a filtering technique to improve the sparsity of the matrices involved in the calculation process of the Paterson-Stockmeyer (PS) scheme. Based on the error analysis considering the truncation error and the error introduced by filtering, the proposed algorithm can obtain similar accuracy as the original PS scheme but is more efficient than it. For the near-sparse matrix power series, the proposed method is also more efficient than the MATLAB built-in codes.

## 1. Introduction

The research on efficient and accurate algorithms for matrix functions has become one of the basic topics in modern scientific computing. Matrix function calculation is involved in computational mathematics, computational mechanics, computational physics, big data analysis, computational finance, and other aspects. Common basic matrix functions include: matrix exponential $e^{\boldsymbol{Z}}$, matrix logarithm $\log \boldsymbol{Z}$, matrix sine / cosine $\sin \boldsymbol{Z} / \cos(\boldsymbol{Z})$, matrix power $\boldsymbol{Z}^p$ ( $p$ is an arbitrary rational number), matrix sign $\operatorname{sign}(\boldsymbol{Z})$, and so on [1]. However, most algorithms proposed for these matrix functions involve the matrix polynomial or the matrix power series [2, 3]. On the one hand, most matrix functions can be expressed as the Taylor series, and it is natural to approximate these matrix functions by using Taylor series. On the other hand, the matrix functions can also be approximated by using the interpolation method, the rational function approximation, or the Chebyshev polynomial approximation, and all these methods involve the matrix polynomial. In a word, the computation of matrix polynomial or matrix power series is a basic issue in the matrix function computation.

There have been many algorithms proposed for matrix polynomial and matrix power series. Among these, the Qin scheme (also called the Horner scheme) and the Paterson-Stockmeyer (PS) scheme are two widely used methods. The PS scheme is an improved version of the Qin scheme and is more efficient for the matrix power series [2, 3]. However, both these methods do not consider the sparsity of the matrix. In many physical and engineering problems, the matrix involved is often highly sparse. Many scholars have found that for these highly sparse matrices, their matrix functions have the property of localization, that is, there are a small number of nonnegligible elements in the matrix functions, while there are a large number of elements very close to zero around them [4, 5]. In this paper, these matrices are called near-sparse matrices. The near-sparse matrix is dense, even full, but most elements in it are very close to zero so dropping out the near-zero elements can generate a sparse approximation to the original matrix.

Since the localization of matrix function was observed, more and more research have been done. Not matrix functions of all the sparse matrices are near-sparse. Under what condition the matrix function is sparse, and how sparse is its approximation have not been fully solved [4-14]. For some specific matrix functions, such as the matrix inversion and the matrix exponential, there are some theoretical estimates about their localization and bandwidths, and these evaluations are often too conservative and difficult to calculate. More accurate and easy to calculate estimates need further research. Readers interested in this issue can read [4] which provided a good review.

The localization of the matrix function can be used for designing more efficient algorithms. In this aspect, the mainstream idea is to improve the computational efficiency of matrix function by using the filtering technique, which removes the elements with very small absolute value in the matrix, and obtains a highly sparse matrix to improve the computation efficiency and reduce the memory. How to determine which elements are small enough to be ignored without affecting the accuracy is an important guarantee for the success of the algorithm. Since the behavior of filtering is equivalent to a kind of artificial rounding error, how to filter the matrices without affecting the accuracy will involve the analysis of filtering error, which is a difficult task. From the existing literature, it can be seen that there are many algorithms based on filtering to improve the computation efficiency for matrix functions such as matrix inverse and matrix exponential. However the research on other functions is quite a few, and the research on matrix series and matrix polynomial, which are two basic matrix functions, has not been reported. It will be more fundamental and important to design an efficient algorithm for computing matrix polynomials and matrix series based on the filtering technique. This work is the content of our research.

The structure of this paper is as follows. In section 2, we briefly review the PS scheme for computing matrix polynomials and matrix series. In section 3, the PS scheme combined with filtering (PSSCF) for matrix polynomials is first presented. According to the error analysis, we propose an adaptive filtering threshold to remove as many elements as possible without losing accuracy. On this basis, the PSSCF is further extended to compute matrix series. In section 4, numerical examples are given to verify the correctness and effectiveness of our algorithm. Finally, the conclusion is given in section 5.

# 2. Paterson and Stockmeyer scheme

## 2.1 Matrix polynomial

The considered matrix polynomial is

$$P_N = \sum_{i=0}^{N-1} a_i \boldsymbol{Z}^i = a_0 \boldsymbol{I} + a_1 \boldsymbol{Z} + a_2 \boldsymbol{Z}^2 + a_3 \boldsymbol{Z}^3 + \ldots + a_N \boldsymbol{Z}^{N-1} \tag{1}$$

where $\boldsymbol{Z}$ is a $n \times n$ matrix, $a_i$ is the coefficients of the polynomial. To evaluate the matrix polynomial, there have been many efficient methods, such as the Qin scheme (also called the Honor scheme) and the Paterson and Stockmeyer (PS) scheme [15]. The PS scheme may be the most efficient one which vitalizes the following equation:

$$\begin{aligned} P_N &= \sum_{i=0}^{N-1} a_i \boldsymbol{Z}^i = a_0 \boldsymbol{I} + a_1 \boldsymbol{Z} + a_2 \boldsymbol{Z}^2 + a_3 \boldsymbol{Z}^3 + \ldots + a_N \boldsymbol{Z}^{N-1} \\ &= \sum_{i=0}^{b-1} \sum_{j=0}^{q-1} a_{iq+j} \boldsymbol{Z}^{iq+j} = \sum_{i=0}^{b-1} \left(\boldsymbol{Z}^q\right)^i \sum_{j=0}^{q-1} a_{iq+j} \boldsymbol{Z}^j \\ &= \sum_{i=0}^{b-1} \left(\boldsymbol{Z}^q\right)^i \boldsymbol{B}_i \end{aligned} \tag{2}$$

where

$$\boldsymbol{B}_i = \sum_{j=0}^{q-1} a_{iq+j} \boldsymbol{Z}^j . \tag{3}$$

According to Eqs. (2) and (3), the calculation process of the PS scheme includes three steps: 1) evaluating $\boldsymbol{I}, \boldsymbol{Z}, \boldsymbol{Z}^2, \cdots, \boldsymbol{Z}^{q-1}, \boldsymbol{Z}^q$ ; 2) evaluating $\boldsymbol{B}_i$ in terms of Eq. (3); 3) evaluating $\sum_{i=0}^{b-1} \left(\boldsymbol{Z}^q\right)^i \boldsymbol{B}_i$ . The first step requires $q-1$ matrix multiplications, and the second step does not require matrix multiplications. The third step needs to evaluate a matrix polynomial which can be done by using the Qin scheme as:

$$\boldsymbol{S}_0 = \boldsymbol{B}_{b-1},\quad \boldsymbol{S}_i = \boldsymbol{B}_{b-1-i} + \boldsymbol{Z}_{(q)}\boldsymbol{S}_{i-1},\quad \boldsymbol{P}_N = \boldsymbol{S}_{b-1} = \boldsymbol{B}_0 + \boldsymbol{Z}_{(q)}\boldsymbol{S}_{b-2},\quad \boldsymbol{Z}_{(q)} = \boldsymbol{Z}^q, \tag{4}$$

which require $b-1$ matrix multiplications. Therefore, the PS scheme requires $q+b-2$ matrix multiplications to evaluate the matrix polynomial. If $\boldsymbol{B}_{b-1} = \sum_{j=0}^{q-1} a_{(b-1)q+j}\boldsymbol{Z}^j = a_{(b-1)q}\boldsymbol{I}$, evaluating $\boldsymbol{S}_1$ does not need matrix multiplications, and the number of matrix multiplications required reduces to $q+b-3$. According to [2], in the PS scheme, the minimum number of matrix multiplications, denoted by $\mathrm{mt}(N)$, is

$$\mathrm{mt}(N) = \begin{cases} 0, & N \le 2 \\ q+b-2-g(q,N), & \text{otherwise} \end{cases} \tag{5}$$

where

$$q = \left\lfloor \sqrt{N-1} \right\rfloor,\quad b = \left\lceil \frac{N}{q} \right\rceil,\quad g(q,N) = \begin{cases} 1, & \text{if} \quad N - q(b-1) = 1 \\ 0, & \text{otherwise} \end{cases}. \tag{6}$$

## 2.2 Matrix power series

The considered matrix power series is

$$f(\boldsymbol{Z}) = \sum_{i=0}^{\infty} a_i \boldsymbol{Z}^i . \tag{7}$$

There are many methods proposed to evaluate the matrix power series, such as the Pade approximation and the Chebyshev approximation. One simple way of these method is taking the first $N$ terms to approximate Eq. (7), i.e.,

$$f(\boldsymbol{Z}) \approx P_N = \sum_{i=0}^{N-1} a_i \boldsymbol{Z}^i = a_0\boldsymbol{I} + a_1\boldsymbol{Z} + a_2\boldsymbol{Z}^2 + a_3\boldsymbol{Z}^3 + \dots + a_N\boldsymbol{Z}^{N-1}, \tag{8}$$

which becomes a matrix polynomial. How to determine the parameter $N$ is a key issue. If let $z_i$ be a upper bound to $\left\|\boldsymbol{Z}^i\right\|$, the upper bound to the error of Eq. (8) can be expressed as:

$$\left\| f(\boldsymbol{Z}) - P_N(\boldsymbol{Z}) \right\| \le \sum_{i=N}^{\infty} |a_i| \left\|\boldsymbol{Z}^i\right\| \le \sum_{i=N}^{\infty} |a_i| z_i . \tag{9}$$

For a given error tolerance $\varepsilon_{\mathrm{tol}}$, $N$ can be determined by $\sum_{i=N}^{\infty} |a_i| z_i \le \varepsilon_{\mathrm{tol}}$.

Equation (9) involves the upper bound to $\left\|\boldsymbol{Z}^i\right\|$ which can be evaluated by a simple method,

i.e., $\|Z^i\| \le z_i = \|Z\|^i$. Of course, this type of evaluation is too conservative, and we here introduce Higham's code **normest1**[16, 17] to evaluate it when $i$ is not too large.

The PS scheme reviewed above requires a great number of matrix multiplications. According to [14], the bandwidth of the matrix power of a sparse matrix increases linearly as the power. Hence the matrix polynomial will be quite dense though the matrix $Z$ is highly sparse. However, much previous research found that there are many elements whose absolute values are very small for many common matrix functions [4, 7]. It is natural to filter out those near-zero elements to obtain a sparse approximation and to improve the computational efficiency. In the next section, we will discuss how to filter out the near-zero elements in the computation of the matrix polynomials.

# 3 PS scheme combined with filtering

This section proposed a PS scheme combined with the filtering technique (PSSCF) to evaluate the matrix power series. As any matrix power series can be approximated by a matrix polynomial with $N$ terms, as shown as Eq. (2), the evaluation of the matrix polynomial will be studied first.

## 3.1 Filtering for matrix polynomial

### 3.1.1 Filtering for a matrix

The near-sparse matrix is a type of matrix which is dense but contains a large number of near-zero elements. In the calculation of matrix functions, we often encounter near-sparse matrices. For a near-sparse matrix $C$, we can drop out the near-zero elements to obtain a sparse $\hat{C}$, and in such a way to improve the computational efficiency. A simple filtering method is to give a filtering threshold $\varepsilon_{\mathrm{f}}$, and then the elements with absolute values smaller than $\varepsilon_{\mathrm{f}}$ will be moved from $C$. However, in the computation of matrix functions, it is often needed to give some conditions to the norm of the matrix $\Delta = C - \hat{C}$ where $\Delta$ consists of the dropped elements. If the matrix size is too large, $\|\Delta\|$ may be also very large even $\varepsilon_{\mathrm{f}}$ is small, and in this case, we cannot state that $\hat{C}$ is an acceptable approximation to $C$. In [14], an adaptive filtering method was proposed. This method can adaptively select the filtering threshold to make sure that $\|\Delta\|$ be smaller than any given filtering-norm-threshold (FNT) $\varepsilon_g$. In this paper, the adaptive filtering method will be used for the calculation of matrix power series. The calculation process of this method is shown in Algorithm 1.

**Algorithm 1** This algorithm filters out the near-zero elements in the matrix $\boldsymbol{C}$ to generate a sparse matrix $\hat{\boldsymbol{C}}$ satisfying $\left\|\boldsymbol{C}-\hat{\boldsymbol{C}}\right\| \le \varepsilon_g\left(1+e_{\mathrm{r}}\right)$. $\varepsilon_g > 0$ is the filtering-norm-threshold, and $e_{\mathrm{r}} > 0$ is a given error tolerance.

1: Set $m_0 = 1$, $c_0 = 1$, $i = 0$, and $\boldsymbol{c}_0 = \boldsymbol{C}$;

2: **while** $c_i > \varepsilon_g\left(1+e_{\mathrm{r}}\right)$;

$\varepsilon_{\mathrm{f}}^{(i+1)} = \varepsilon_g / m_i$;

$p = \mathrm{find}\left(\left|\boldsymbol{c}_i\right| > \varepsilon_{\mathrm{f}}^{(i+1)}\right)$, where $p$ is a set of linear indices corresponding to the nonzero entries in $\boldsymbol{c}_i$ whose absolute values are larger than $\varepsilon_{\mathrm{f}}^{(i+1)}$;

$\boldsymbol{c}_i\left(p\right) = 0, \quad \boldsymbol{c}_{i+1} = \boldsymbol{c}_i$;

$c_{i+1} = \left\|\boldsymbol{c}_{i+1}\right\|$;

$m_{i+1} = c_{i+1} / \varepsilon_{\mathrm{f}}^{(i+1)}$;

$i = i + 1$;

**end**

3: $\hat{\boldsymbol{C}} = \boldsymbol{C} - \boldsymbol{c}_i$.

### 3.1.2 Formula of PS scheme with filtering

For the convenience of analysis, let $\boldsymbol{Z}_{(i)} := \boldsymbol{Z}^i$. The first step of the PS scheme can be shown as

$$\boldsymbol{Z}_{(2)} = \boldsymbol{Z}^2,\ \boldsymbol{Z}_{(i+1)} = \boldsymbol{Z}\boldsymbol{Z}_{(i)},\ \boldsymbol{Z}_{(q)} = \boldsymbol{Z}\boldsymbol{Z}_{(q-1)}, \tag{10}$$

which involves many matrix multiplications and the matrix $\boldsymbol{Z}_{(i)}$ becomes denser and denser. If we drop out the near-zero elements in $\boldsymbol{Z}_{(i)}$, the calculation process can be shown as

$$\hat{\boldsymbol{Z}}_{(1)} = \boldsymbol{Z},\ \hat{\boldsymbol{Z}}_{(2)} = \boldsymbol{Z}\hat{\boldsymbol{Z}}_{(1)} - \boldsymbol{\zeta}_2,\ \hat{\boldsymbol{Z}}_{(i+1)} = \boldsymbol{Z}\hat{\boldsymbol{Z}}_{(i)} - \boldsymbol{\zeta}_{i+1} \tag{11}$$

where $\hat{\boldsymbol{Z}}_{(i+1)}$ represents the sparse matrix obtained by filtering $\boldsymbol{Z}\hat{\boldsymbol{Z}}_{(i)}$, and $\boldsymbol{\zeta}_{i+1}$ is the matrix consisting of the dropped-out elements. Combining $\boldsymbol{Z}_{(i+1)} = \hat{\boldsymbol{Z}}_{(i+1)} + \boldsymbol{\zeta}_{i+1}$ with Eq. (11) yields

$$\hat{\boldsymbol{Z}}_{(i)} = \boldsymbol{Z}^i - \sum_{k=2}^{i} \boldsymbol{Z}^{i-k}\boldsymbol{\zeta}_k, \tag{12}$$

where $\Delta\boldsymbol{Z}_{(i)} := \sum_{k=2}^{i} \boldsymbol{Z}^{i-k}\boldsymbol{\zeta}_k$ represents the error between $\hat{\boldsymbol{Z}}_{(i)}$ and $\boldsymbol{Z}^i$. Substituting Eq. (12) into Eq. (3), we have

$$\begin{aligned} \boldsymbol{B}_i &= \sum_{j=0}^{q-1} a_{iq+j} \boldsymbol{Z}^j = \sum_{j=0}^{q-1} a_{iq+j} \hat{\boldsymbol{Z}}_{(j)} + \sum_{j=0}^{q-1} \left( a_{iq+j} \sum_{k=2}^{j} \boldsymbol{Z}^{j-k} \boldsymbol{\zeta}_k \right) \\ &= \hat{\boldsymbol{B}}_i + \sum_{j=2}^{q-1} \left( \sum_{k=j}^{q-1} a_{iq+k} \boldsymbol{Z}^{k-j} \right) \boldsymbol{\zeta}_j = \hat{\boldsymbol{B}}_i + \Delta \boldsymbol{B}_i \end{aligned} . \tag{13}$$

where $\Delta\boldsymbol{B}_i$ is the error between $\hat{\boldsymbol{B}}_i$ and $\boldsymbol{B}_i$. As the calculation of $\boldsymbol{B}_i$ does not require matrix multiplications, filtering is also needed. But the third step of the PS scheme requires matrix multiplications, and hence also needs filtering. The third step is calculated by using the Qin scheme. If the filtering is considered, the Qin scheme should be modified as

$$\hat{\boldsymbol{S}}_0 = \hat{\boldsymbol{B}}_{b-1}, \quad \hat{\boldsymbol{S}}_i = \hat{\boldsymbol{B}}_{b-1-i} + \hat{\boldsymbol{Z}}_{(q)} \hat{\boldsymbol{S}}_{i-1} - \boldsymbol{\Sigma}_i, \quad \hat{P}_N = \hat{\boldsymbol{S}}_{b-1} = \hat{\boldsymbol{B}}_0 + \hat{\boldsymbol{Z}}_{(q)} \hat{\boldsymbol{S}}_{b-2} - \boldsymbol{\Sigma}_{b-1} \tag{14}$$

where $\hat{\boldsymbol{S}}_i$ represents the sparse matrix obtained by filtering $\hat{\boldsymbol{B}}_{b-1-i} + \hat{\boldsymbol{Z}}_{(q)} \hat{\boldsymbol{S}}_{i-1}\ \boldsymbol{Z}_{(i+1)}$, $\boldsymbol{\Sigma}_i$ consists of the dropped out elements, and $\hat{P}_N$ is the approximation to $P_N$ . When $\boldsymbol{B}_{b-1} = a_{(b-1)q} \boldsymbol{I}$ , $\hat{\boldsymbol{B}}_{b-2} + \hat{\boldsymbol{Z}}_{(q)} \hat{\boldsymbol{S}}_0 = \hat{\boldsymbol{B}}_{b-2} + \hat{\boldsymbol{Z}}_{(q)} a_{(b-1)q}$ does not require the matrix multiplication, and hence also does not need filtering, i.e., $\boldsymbol{\Sigma}_{b-1} = \boldsymbol{0}$ .

### 3.1.3 Error introduced by filtering

According to Eq. (14), we have

$$\hat{\boldsymbol{S}}_i = \hat{\boldsymbol{B}}_{b-1-i} + \hat{\boldsymbol{Z}}_{(q)} \hat{\boldsymbol{S}}_{i-1} - \boldsymbol{\Sigma}_i = \boldsymbol{B}_{b-1-i} - \Delta\boldsymbol{B}_{b-1-i} + \left( \boldsymbol{Z}^q - \Delta\boldsymbol{Z}_{(q)} \right) \left( \boldsymbol{S}_{i-1} - \Delta\boldsymbol{S}_{i-1} \right) - \boldsymbol{\Sigma}_i = \boldsymbol{S}_i - \Delta\boldsymbol{S}_i . \tag{15}$$

To simplify the analysis, we only consider the case that the errors $\left\| \Delta\boldsymbol{Z}_{(q)} \right\|$ and $\left\| \Delta\boldsymbol{S}_{i-1} \right\|$ are so small that $\Delta\boldsymbol{Z}_{(q)} \Delta\boldsymbol{S}_{i-1}$ can be ignored. In this case, $\Delta\boldsymbol{S}_i$ can be written as

$$\Delta\boldsymbol{S}_i = \boldsymbol{Z}_{(q)} \Delta\boldsymbol{S}_{i-1} + \boldsymbol{\Phi}_i, \quad \boldsymbol{\Phi}_i = \Delta\boldsymbol{B}_{b-1-i} + \Delta\boldsymbol{Z}_{(q)} \boldsymbol{S}_{i-1} + \boldsymbol{\Sigma}_i \tag{16}$$

and

$$\Delta\boldsymbol{S}_0 = \boldsymbol{S}_0 - \hat{\boldsymbol{S}}_0 = \Delta\boldsymbol{B}_{b-1} . \tag{17}$$

Combining Eq. (16) with Eq. (17), $\Delta\boldsymbol{S}_i$ can also be rewritten as

$$\Delta\boldsymbol{S}_i = \sum_{j=0}^{i} \boldsymbol{Z}_{(q)}^{i-j} \boldsymbol{\Phi}_j, \quad \boldsymbol{\Phi}_0 = \Delta\boldsymbol{S}_0 = \Delta\boldsymbol{B}_{b-1} . \tag{18}$$

From (18), and noting that $\left\| \boldsymbol{Z}_{(q)} \right\| = \left\| \boldsymbol{Z}^q \right\|$, we have

$$\left\| \Delta\boldsymbol{S}_{b-1} \right\| \le \sum_{j=0}^{b-1} \left\| \boldsymbol{Z}^{q(b-1-j)} \right\| \left\| \boldsymbol{\Phi}_j \right\| \tag{19}$$

where

$$\left\|\boldsymbol{\Phi}_j\right\| \le \left\|\Delta\boldsymbol{B}_{b-1-j}\right\| + \left\|\Delta\boldsymbol{Z}_{(q)}\right\|\left\|\boldsymbol{S}_{j-1}\right\| + \left\|\boldsymbol{\Sigma}_j\right\| \tag{20}$$

From (13),

$$\left\|\Delta\boldsymbol{B}_i\right\| \le \sum_{j=2}^{q-1}\beta_{i,j}\left\|\boldsymbol{\zeta}_j\right\|, \quad \beta_{i,j} = \sum_{k=j}^{q-1}\left|a_{iq+k}\right|\left\|\boldsymbol{Z}^{k-j}\right\| \tag{21}$$

From (12),

$$\left\|\Delta\boldsymbol{Z}_{(q)}\right\| \le \sum_{k=2}^{q}\left\|\boldsymbol{Z}^{q-k}\right\|\left\|\boldsymbol{\zeta}_k\right\| \tag{22}$$

Combining Eqs. (19)-(22) yields

$$\left\|\Delta\boldsymbol{S}_{b-1}\right\| \le \sum_{k=2}^{q} d_k\left\|\boldsymbol{\zeta}_k\right\| + \sum_{j=1}^{b-1-g(q,N)}\left\|\boldsymbol{Z}^{q(b-1-j)}\right\|\left\|\boldsymbol{\Sigma}_j\right\| \tag{23}$$

where $g(q,N)$ is shown as Eq. (6), and

$$d_k = \sum_{j=0}^{b-1}\left\|\boldsymbol{Z}^{q(b-1-j)}\right\|\left(\beta_{b-1-j,k} + \left\|\boldsymbol{Z}^{q-k}\right\|\left\|\boldsymbol{S}_{j-1}\right\|\right), \ \left\|\boldsymbol{S}_{-1}\right\| = \beta_{b-1-j,q} = 0\,. \tag{24}$$

From Eqs. (3) and (4), we have

$$\left\|\boldsymbol{B}_i\right\| = \left\|\sum_{j=0}^{q-1} a_{iq+j}\boldsymbol{Z}^j\right\| \le \sum_{j=0}^{q-1}\left|a_{iq+j}\right|\left\|\boldsymbol{Z}^j\right\|, \tag{25}$$

and

$$\left\|\boldsymbol{S}_0\right\| = \left\|\boldsymbol{B}_{b-1}\right\|, \quad \left\|\boldsymbol{S}_i\right\| = \left\|\boldsymbol{B}_{b-1-i} + \boldsymbol{Z}_{(q)}\boldsymbol{S}_{i-1}\right\| \le \left\|\boldsymbol{B}_{b-1-i}\right\| + \left\|\boldsymbol{Z}^q\right\|\left\|\boldsymbol{S}_{i-1}\right\| \tag{26}$$

Using Eqs. (25)-(26) and Eq. (24), we can evaluate $d_k$ .

### 3.1.4 Adaptive FNT

Equation (23) is the error bound to the PS scheme combined with filtering. For a given error tolerance $\varepsilon$, if we require that

$$\left\|\Delta\boldsymbol{S}_{b-1}\right\| \le \sum_{k=2}^{q} d_k\left\|\boldsymbol{\zeta}_k\right\| + \sum_{j=1}^{b-1-g(q,N)}\left\|\boldsymbol{Z}^{q(b-1-j)}\right\|\left\|\boldsymbol{\Sigma}_j\right\| < \varepsilon \tag{27}$$

the upper bounds to $\left\|\boldsymbol{\zeta}_k\right\|$ and $\left\|\boldsymbol{\Sigma}_j\right\|$ should be constrained.

In Eq. (23), $d_k\left\|\boldsymbol{\zeta}_k\right\|$ represents the error introduced by the filtering matrix $\boldsymbol{\zeta}_k$ , and $\left\|\boldsymbol{Z}^{q(b-1-j)}\right\|\left\|\boldsymbol{\Sigma}_j\right\|$ the error introduced by the filtering matrix $\boldsymbol{\Sigma}_j$ . If we require that the influences of different filtering matrices on the last error $\Delta\boldsymbol{S}_{b-1}$ are the same, we can derive the filtering-norm-threshold for $\boldsymbol{\zeta}_k$ and $\boldsymbol{\Sigma}_j$ that

$$\left\|\boldsymbol{\zeta}_k\right\| \le \zeta_{k,\mathrm{p}} = \frac{\varepsilon}{d_k\left(q+b-2-g(q,N)\right)}, \quad \left\|\boldsymbol{\Sigma}_j\right\| \le \sigma_{j,\mathrm{p}} = \frac{\varepsilon}{\left\|\boldsymbol{Z}^{q(b-1-j)}\right\|\left(q+b-2-g(q,N)\right)} \tag{28}$$

where $d_k$ can be evaluated in terms of Eq. (24).

## 3.2 Filtering for matrix power series

The matrix power series $f\left(\boldsymbol{Z}\right)=\sum_{i=0}^{\infty}a_i\boldsymbol{Z}^i$ can be approximated by the sum of the first $N$ terms, i.e., $P_N\left(\boldsymbol{Z}\right)=\sum_{i=0}^{N-1}a_i\boldsymbol{Z}^i$, which is a matrix polynomial and can be computed by using the PSSCF. However, the error introduced by truncating the series should be considered to determine the filtering-norm-thresholds.

### 3.2.1 Filtering-norm-thresholds considering the truncation error

For the matrix power series, the error can be written as

$$f\left(\boldsymbol{Z}\right)-\hat{\boldsymbol{S}}_{b-1}=f\left(\boldsymbol{Z}\right)-P_N\left(\boldsymbol{Z}\right)+P_N\left(\boldsymbol{Z}\right)-\hat{\boldsymbol{S}}_{b-1}=f\left(\boldsymbol{Z}\right)-P_N\left(\boldsymbol{Z}\right)+\Delta\boldsymbol{S}_{b-1}, \tag{29}$$

for which we have

$$\left\|f\left(\boldsymbol{Z}\right)-\hat{\boldsymbol{S}}_{b-1}\right\|\le\left\|f\left(\boldsymbol{Z}\right)-P_N\left(\boldsymbol{Z}\right)\right\|+\left\|\Delta\boldsymbol{S}_{b-1}\right\|. \tag{30}$$

For a given error tolerance $\varepsilon_{\text{tol}}$, the error should satisfy that

$$\left\|\Delta\boldsymbol{S}_{b-1}\right\|+\left\|f\left(\boldsymbol{Z}\right)-P_N\left(\boldsymbol{Z}\right)\right\|\le\sum_{k=2}^{q}d_k\left\|\boldsymbol{\zeta}_k\right\|+\sum_{j=1}^{b-1-g(q,N)}\left\|\boldsymbol{Z}^{q(b-1-j)}\right\|\left\|\boldsymbol{\Sigma}_j\right\|+\left\|f\left(\boldsymbol{Z}\right)-P_N\left(\boldsymbol{Z}\right)\right\|<\varepsilon_{\text{tol}}, \tag{31}$$

where, $\sum_{k=2}^{q}d_k\left\|\boldsymbol{\zeta}_k\right\|+\sum_{j=1}^{b-1-g(q,N)}\left\|\boldsymbol{Z}^{q(b-1-j)}\right\|\left\|\boldsymbol{\Sigma}_j\right\|$ is the error introduced by filtering, and $\left\|f\left(\boldsymbol{Z}\right)-P_N\left(\boldsymbol{Z}\right)\right\|$ is the truncation error. If the upper bounds to the errors introduced by filtering and truncating are required to be the same, we have

$$\left\|f\left(\boldsymbol{Z}\right)-P_N\left(\boldsymbol{Z}\right)\right\|\le0.5\varepsilon_{\text{tol}}. \tag{32}$$

What should be noted is the truncation error is controlled by $N$ which is unknown before calculating. Combining Eqs. (9) and (32), we have

$$\left\|f\left(\boldsymbol{Z}\right)-P_N\left(\boldsymbol{Z}\right)\right\|\le\sum_{i=N}^{\infty}\left|a_i\right|\left\|\boldsymbol{Z}^i\right\|\le\sum_{i=N}^{\infty}\left|a_i\right|z_i\le0.5\varepsilon_{\text{tol}}, \tag{33}$$

which yields

$$N=\Theta\left(\varepsilon_{\text{tol}},\boldsymbol{a},\boldsymbol{z}\right):=\min\left(m:\ \sum_{i=m}^{\infty}\left|a_i\right|z_i\le0.5\varepsilon_{\text{tol}},\ \text{and,}\ m\in\mathbb{N}\right) \tag{34}$$

where $\boldsymbol{z}=\left(z_0,z_1,\cdots,z_i,\cdots\right)$, and $\boldsymbol{a}=\left(a_0,a_1,\cdots,a_i,\cdots\right)$ is an infinite vector consisting of the coefficients of the matrix power series. According to Eq. (34), $N$ can be evaluated. Combining (31) with (33), we have

$$\sum_{k=2}^{q} d_k \|\boldsymbol{\zeta}_k\| + \sum_{j=1}^{b-1-g(q,N)} \left\|\boldsymbol{Z}^{q(b-1-j)}\right\| \|\boldsymbol{\Sigma}_j\| \le \varepsilon_{\text{tol}} - \sum_{i=N}^{\infty} |a_i| z_i \ . \tag{35}$$

Again, requiring the influences of the filtering matrices are the same yields the FNT that

$$\|\boldsymbol{\zeta}_k\| \le \zeta_{k,s} = \frac{\varepsilon_{\text{tol}} - \sum_{i=N}^{\infty} |a_i| z_i}{d_k \left(q + b - 2 - g(q,N)\right)}, \quad \|\boldsymbol{\Sigma}_j\| \le \sigma_{j,s} = \frac{\varepsilon_{\text{tol}} - \sum_{i=N}^{\infty} |a_i| z_i}{\left\|\boldsymbol{Z}^{q(b-1-j)}\right\| \left(q + b - 2 - g(q,N)\right)} \tag{36}$$

### 3.2.2 Evaluation of the norm of matrix power

It can be seen from the above analysis, the norm $\|\boldsymbol{Z}^i\|$ has a great effect on the FNT. In [17], Higham proposed an efficient method and the corresponding MATLAB code **normest1** for the evaluation of $\|\boldsymbol{Z}^i\|$. After numerical simulations, it was observed that the computation time of $\|\boldsymbol{Z}^i\|$ using the **normest1** increases as $i$. In this paper, if $i \le 10$, we use **normest1** to evaluate $\|\boldsymbol{Z}^i\|$; else, we use the following lemma [18].

**Lemma 1** Let $\|\boldsymbol{Z}^i\|$ be known for $i = 0, 1, \cdots, K$, and $K \ge 2$ be an even number, then for any $i \ge \left\lceil \frac{K+1}{2} \right\rceil$, we have $\|\boldsymbol{Z}^i\| \le \alpha_K^i$, where $\alpha_K = \max\left\{ \|\boldsymbol{Z}^k\|^{\frac{1}{k}},\ k = \left\lceil \frac{K}{2} \right\rceil, \cdots, K \right\}$.

Lemma 1 is a simple corollar of the Theorem 2 in [18]. Using Higham's method and Lemma 1, the evaluation of $z_i$ can be shown as

$$z_i = \begin{cases} \text{normest1}\left(\boldsymbol{Z}^i\right), & i \le 10 \\ \alpha_{10}^i, \ i > 10 \end{cases}, \quad \alpha_{10} = \max\left\{ \|\boldsymbol{Z}^k\|^{\frac{1}{k}},\ k = 5:10 \right\} \tag{37}$$

### 3.2.3 PSSCF for matrix power series

**Algorithm 2** This algorithm evaluates the approximation $\hat{f}(\boldsymbol{Z})$ to the matrix power series $f(\boldsymbol{Z}) = \sum_{i=0}^{\infty} a_i \boldsymbol{Z}^i$ satisfying $\left\| f(\boldsymbol{Z}) - \hat{f}(\boldsymbol{Z}) \right\| \le \varepsilon_{\text{tol}}$, where $\varepsilon_{\text{tol}} > 0$ is a given error tolerance.

1: Calculate $\boldsymbol{z} = (z_0, z_1, \cdots, z_i, \cdots)$ with formula (37);

2: $N = \Theta(\varepsilon_{\text{tol}}, \boldsymbol{a}, \boldsymbol{z})$;

3: $q = \left\lfloor \sqrt{N-1} \right\rfloor,\ b = \left\lceil \frac{N}{q} \right\rceil$;

4: **for** $i = 1:(q-1)$;

5: $\tilde{\boldsymbol{Z}}_{(i+1)} = \boldsymbol{Z}\hat{\boldsymbol{Z}}_{(i)}$;

6: Estimate $\zeta_{i+1,\mathrm{s}}$ according to (36), filter out $\tilde{\boldsymbol{Z}}_{(i+1)}$ with **Algorithm 1** to get $\hat{\boldsymbol{Z}}_{(i+1)}$, satisfying $\left\|\tilde{\boldsymbol{Z}}_{(i+1)} - \hat{\boldsymbol{Z}}_{(i+1)}\right\| \le \zeta_{i+1,\mathrm{s}}$;

7: **end**

8: $\hat{\boldsymbol{S}}_0 = \hat{\boldsymbol{B}}_{b-1} = \sum_{j=0}^{q-1} a_{(b-1)q+j}\hat{\boldsymbol{Z}}_{(j)}$;

9: **for** $i = 1:(q-1)$;

10: $\hat{\boldsymbol{B}}_{b-1-i} = \sum_{j=0}^{q-1} a_{(b-1-i)q+j}\hat{\boldsymbol{Z}}_{(j)}$;

11: $\tilde{\boldsymbol{S}}_i = \hat{\boldsymbol{B}}_{b-1-i} + \hat{\boldsymbol{Z}}_{(q)}\hat{\boldsymbol{S}}_{i-1}$;

12: Estimate $\sigma_{i,\mathrm{s}}$ according to (36), filter out $\tilde{\boldsymbol{S}}_i$ with **Algorithm 1** to get $\hat{\boldsymbol{S}}_i$, satisfying $\left\|\tilde{\boldsymbol{S}}_i - \hat{\boldsymbol{S}}_i\right\| \le \sigma_{i,\mathrm{s}}$;

13: **end**

14: $\hat{f}(\boldsymbol{Z}) = \hat{\boldsymbol{S}}_{b-1}$

If the matrix polynomial is evaluated, the second step in Algorithm 2 should be ignored due to the number $N$ is known, and the FNTs of $\zeta_{i+1,\mathrm{s}}$ and $\sigma_{i,\mathrm{s}}$ should be replaced with $\zeta_{i+1,\mathrm{p}}$ and $\sigma_{i,\mathrm{p}}$, as shown in (28).

# 4. Numerical examples

Two types of matrix functions, i.e., the matrix exponential and the matrix cosine, are used to test the effectiveness and efficiency of the proposed method. These two matrix functions can be expressed as the Taylor series which converge for any matrices. For the matrix exponential, the original PS scheme without filtering, the proposed PSSCF, and the MATLAB code **expm** are used; and for the matrix exponential, the original PS scheme without filtering, the proposed PSSCF, and the MATLAB code **cosm** are used. Both the original PS scheme and the proposed PSSCF have been implemented by using MATLAB, and the corresponding MATLAB codes have been uploaded to https://www.rocewea.com/10.html. We performed the experiments by using the computer with Microsoft Windows 11 22H2, AMD Ryzen 7 5800H with Radeon Graphics @3.20 GHz, and 15.9GB of RAM. And the MATLAB version is MATLAB R2021b.

We used the following sets of matrices for testing:

1. Ten 10000×10000 adjacency matrices $\boldsymbol{B}_i$ were downloaded from https://networkrepository.com. These matrices were modified to be $\boldsymbol{A}_i = \boldsymbol{I} - 0.5\boldsymbol{B}_i / \rho(\boldsymbol{B}_i)$, which were used as the test matrices. The matrices $\boldsymbol{A}_i$ can also be downloaded from https://www.rocewea.com/10.html.
2. Random banded matrices 1000×1000 with different bandwidths were used. The lower and upper bandwidths of these matrices are the same, and denoted by $b$, and $b$ ranges from 1 to 15. The random numbers in these matrices obey the standard normal distribution. For each bandwidth, 30 random matrices are generated.

## 4.1 Matrix exponential

For the sparse matrix $\boldsymbol{A}_i$, the matrix exponential can be written as

$$e^{\boldsymbol{A}_i} = \sum_{j=0}^{\infty} \frac{1}{j!} \boldsymbol{A}_i^j . \tag{38}$$

To test the efficiency and accuracy of the proposed PSSCF, the original PS scheme, the proposed PSSCF, and the MATLAB built-in code **expm** are used to evaluate the matrix exponential. The **expm** result is seen as the benchmark solution. The errors of the original PS and proposed methods are evaluated by

$$er_i = \frac{\left\| \boldsymbol{X}_i - \mathrm{expm}(\boldsymbol{A}_i) \right\|}{\left\| \mathrm{expm}(\boldsymbol{A}_i) \right\|} \tag{39}$$

where $\boldsymbol{X}_i$ represents the result evaluated with the original PS or the proposed method. For both the original PS scheme and the proposed PSSCF, the error tolerances are set to be $\varepsilon_{\mathrm{tol}} = 10^{-14}$.

Table 1 Comparison between different algorithms used to compute the matrix exponential experiments of the first set of test matrices.

| $\boldsymbol{A}_i$ | Computational time (s) | | | $er_i$ | |
|---|---|---|---|---|---|
| | PSSCF | PS | **expm** | PSSCF | PS |
| 1 | 0.06 | 0.06 | 49.81 | 5.33E-13 | 5.33E-13 |
| 2 | 1.60 | 3.75 | 50.67 | 2.72E-13 | 2.72E-13 |
| 3 | 3.88 | 12.31 | 20.11 | 2.70E-15 | 2.70E-15 |
| 4 | 0.05 | 0.07 | 11.06 | 3.90E-13 | 3.90E-13 |
| 5 | 0.05 | 0.06 | 49.65 | 5.50E-13 | 5.50E-13 |
| 6 | 0.23 | 0.30 | 49.04 | 8.74E-13 | 8.74E-13 |
| 7 | 0.07 | 0.08 | 50.25 | 4.88E-13 | 4.88E-13 |
| 8 | 0.09 | 0.14 | 50.45 | 4.83E-13 | 4.83E-13 |
| 9 | 0.10 | 0.13 | 49.30 | 4.91E-13 | 4.91E-13 |
| 10 | 0.14 | 0.22 | 50.01 | 4.85E-13 | 4.85E-13 |

The computational times and relative errors of different algorithms used for the first set of test matrices are compared in Table 1. As shown, both the original PS scheme and the proposed PSSCF can evaluate the considered test matrices precisely and efficiently. The accuracy of PSSCF is equivalent to that of PS, but the efficiency of PSSCF is better than that of PS, which is particularly obvious in the third test matrix $\boldsymbol{A}_3$, where the computational time of PSSCF is only one-fourth of that of PS.

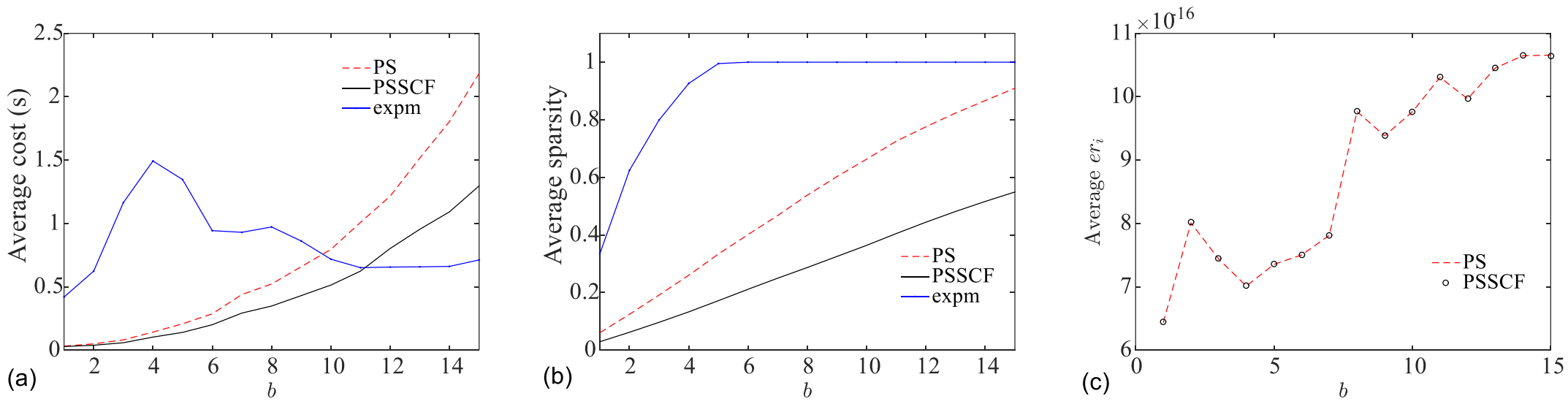

Figure 1. Comparison of different algorithms in terms of the average computational times (a), the average sparsity of the computed matrix exponentials (b), and the relative errors (c).

Figures 1 (a) compares the average computational times that different algorithms take to compute the 30 random banded matrices with each bandwidth. The approximate matrices obtained by different methods have different sparsity, i.e., the ratio of the number of non-zero elements to the square of the dimension of the matrix. The average sparsities of the 30 results obtained by different algorithms are compared in Fig. 1(b). The average errors of PS and PSSCF are compared in Fig. 1(c). As shown in Fig. 1, both the PS and the PSSCF can produce results with high accuracy. The sparsity of the matrix exponential approximated by PSSCF is smaller than that by PS or **expm**. PSSCF is always more efficient than PS. Compared with **expm**, PSSCF performs more efficient when $b \le 11$, and less efficient when $b > 11$, as the matrix exponential is no longer near-sparse when $b > 11$ in terms of Fig. 1(b).

## 4.2 Matrix cosine

For the sparse matrix $\boldsymbol{A}_i$, the matrix exponential can be written as

$$\cos \boldsymbol{A}_i = \sum_{j=0}^{\infty} (-1)^j \frac{\boldsymbol{A}_i^{2j}}{(2j)!}. \tag{40}$$

To test the efficiency and accuracy of the proposed PSSCF, the original PS scheme, the proposed PSSCF, and the MATLAB built-in code **cosm** are used to evaluate the matrix cosine. The **cosm** result is seen as the benchmark solution. The errors of the original PS and proposed methods are evaluated by

$$er_i = \frac{\left\| \boldsymbol{X}_i - \mathrm{cosm}(\boldsymbol{A}_i) \right\|}{\left\| \mathrm{cosm}(\boldsymbol{A}_i) \right\|}, \tag{41}$$

where $\boldsymbol{X}_i$ represents the matrix cosine evaluated with the original PS or the proposed method. For both the original PS scheme and the proposed PSSCF, the error tolerances are set to be $\varepsilon_{\mathrm{tol}} = 10^{-14}$.

Table 2 Comparison between different algorithms used to compute the matrix cosine experiments of the first set of test matrices.

| $\boldsymbol{A}_i$ | Computational time (s) | | | $er_i$ | |
|---|---|---|---|---|---|
| | PSSCF | PS | **cosm** | PSSCF | PS |
| 1 | 0.04 | 0.04 | 30.19 | 1.73E-15 | 3.80E-16 |
| 2 | 1.57 | 4.02 | 77.18 | 1.46E-14 | 1.42E-14 |
| 3 | 4.39 | 15.79 | 121.86 | 1.59E-14 | 1.48E-14 |
| 4 | 0.04 | 0.04 | 28.21 | 2.68E-15 | 5.93E-16 |
| 5 | 0.04 | 0.05 | 30.19 | 1.63E-15 | 6.31E-16 |
| 6 | 0.25 | 0.22 | 36.59 | 3.92E-15 | 2.36E-15 |
| 7 | 0.06 | 0.08 | 31.66 | 2.06E-15 | 3.57E-16 |
| 8 | 0.09 | 0.09 | 33.78 | 2.41E-15 | 4.43E-16 |
| 9 | 0.08 | 0.10 | 33.87 | 2.28E-15 | 6.42E-16 |
| 10 | 0.14 | 0.17 | 36.67 | 2.58E-15 | 5.21E-16 |

The computational times and relative errors of different algorithms are compared in Table 2. As shown, both the original PS scheme and the proposed PSSCF can evaluate the considered test matrices precisely and efficiently. The accuracy of PSSCF is equivalent to that of PS. For the test matrices $\boldsymbol{A}_2$ and $\boldsymbol{A}_3$, the computational times of PSSCF are far smaller than that of PS, which shows the filtering can significantly improve the computational efficiency of PS. For other test matrices, the computational efficiency of PSSCF and PS is equivalent, which may be because the considered matrices are too sparse to reflect the advantages of the filtering. But it also shows that the filtering does not weaken the computational efficiency and accuracy of the original PS.

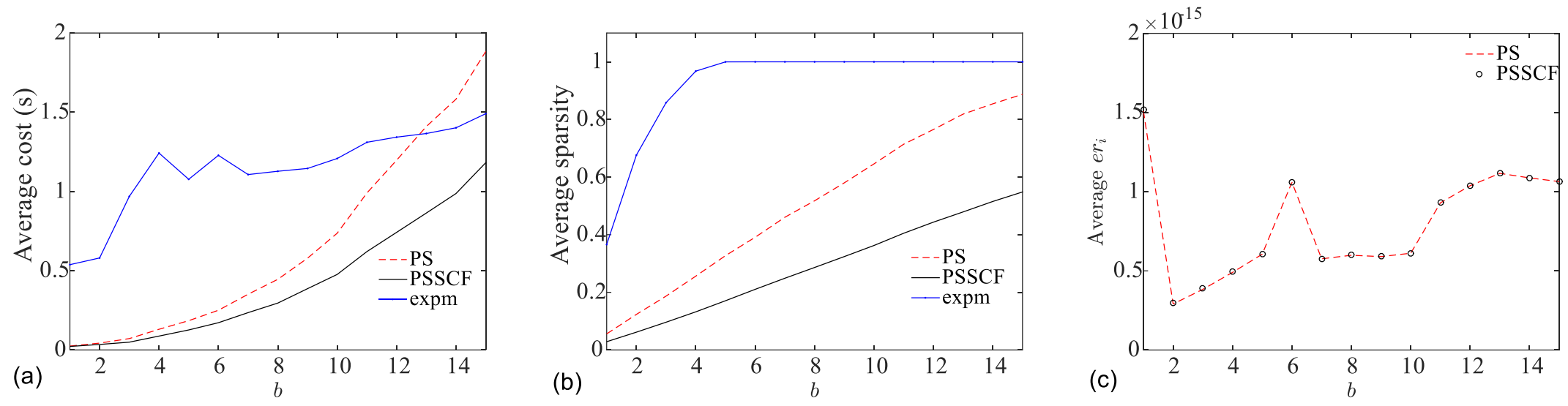

Figure 2. Comparison of different algorithms in terms of the average computational times (a), the average sparsity of the computed matrix exponentials (b), and the relative errors (c).

Figures 2 (a) compares the average computational times that different algorithms take to compute the matrix cosines of the 30 random banded matrices with each bandwidth. The average sparsities of the 30 results obtained by different algorithms are compared in Fig. 2(b). The average errors of PS and PSSCF are compared in Fig. 2(c). As shown in Fig. 2 again, both the PS and the PSSCF can produce results with high accuracy. The sparsity of the matrix cosines approximated by PSSCF is smaller than that by PS or **cosm**. PSSCF is always more efficient than PS or **cosm**.

# 5. Conclusion

A competitive modification of the PS scheme has been proposed by using the filtering technique. The filtering-norm-threshold was based on the error analysis considering the truncation error and the error introduced by filtering, leading to a new algorithm for the matrix power series which is near-sparse. It was verified by numerical experiments on evaluating the matrix exponentials and matrix cosines of two sets of sparse matrices that the proposed method can obtain similar accuracy as the original PS scheme, but is more efficient than the PS scheme. For the near-sparse matrix power series, the proposed method is also more efficient than some state-of-the-art codes, such as **expm** and **cosm**.